\documentclass[11pt,twoside]{amsart}
\usepackage[T1]{fontenc}
\usepackage[utf8]{inputenc}
\usepackage[english]{babel}
\usepackage{graphicx}
\usepackage{amsmath}
\usepackage{amssymb}
\usepackage{enumerate}
\usepackage[all,cmtip]{xy}
\usepackage{mathrsfs}

\usepackage{amsthm}
\usepackage{amsfonts}
\usepackage{changepage}

\newtheorem{thm1}{Theorem}[section]
\newtheorem{theorem}[thm1]{Theorem}
\newtheorem{lemma}[thm1]{Lemma}

\newtheorem{proposition}[thm1]{Proposition}

\newtheorem{mainthm}{Theorem}

\theoremstyle{definition}
\newtheorem{definition}[thm1]{Definition}

\theoremstyle{remark}

\title[On Carleman-Sobolev Classes for exponents $0<p<1$]{Independence of derivatives in Carleman-Sobolev Classes for exponents $0<p<1$}
\author{Aron Wennman}

\renewcommand{\phi}{\varphi}
\renewcommand{\epsilon}{\varepsilon}

\newcommand{\R}{\mathbb{R}}

\renewcommand{\and}{\hspace{8pt} \text{and} \hspace{8pt}}
\newcommand{\subsec}[1]{\par\vspace{\baselineskip}
\noindent\textbf{#1.}}

\newcommand{\dd}{\,d}
\newcommand{\nl}{\left\lVert}
\newcommand{\nr}{\right\rVert}

\newcommand{\M}{\mathcal{M}}
\newcommand{\MS}{{\mathcal{M}_1}}
\newcommand{\product}[1]{\mu_{\mathcal{#1}}}
\newcommand{\norm}[1]{\nl #1\nr_{\M}}
\renewcommand{\S}{\mathcal{S}}

\newcommand{\C}{\mathcal{C}}

\newcommand{\supp}{\operatorname{supp}}
\newcommand{\W}{W^{p}_{\! \M}}
\newcommand{\Winf}{\mathscr{W}^{p}_{\! \M}}
\newcommand{\Winfshift}{\mathscr{W}^{p}_{\! \MS}}

\DeclareMathOperator{\id}{id}
\DeclareMathOperator{\im}{im}

\numberwithin{equation}{section}
\address{Department of Mathematics, KTH Royal Institute of Technology, SE-100 44 Stockholm, Sweden}

\email{aronw@math.kth.se}
\begin{document}

\begin{abstract}
We continue the study of Carleman-Sobolev classes from previous joint work with G.~Behm. 
We consider spaces denoted by $\W$, defined as abstract completions of sets of smooth functions with respect to a weighted Sobolev-flavoured norm involving derivatives of all orders. 
Previously we showed that these classes behave very differently on two sides of a condition on the weight sequence $\M$. Here we prove a conjecture made in the previous work;  under some regularity assumptions on the weight, we show that on one side of the condition there will be a complete independence between derivatives, expressed as
$$
\W\cong L^p\oplus W^p_{\MS}
$$
where $\MS$ is the shifted sequence. On the other side, we already know that one can embed $\W$ into $C^{\infty}(\R)$. Thus this is an instance of a kind of phase transition.
\end{abstract}

\maketitle

\section*{Introduction}

\subsec{Sobolev spaces for $0<p<1$}
In \cite{Peetre}, Jaak Peetre initiates the study of a certain version of Sobolev spaces for small $L^p$-exponents, i.e.~$p$ in the range $0<p<1$. 
He considers the abstract completion of $C_0^{\infty}(\R)$ with respect to the quasi-norm
$$
\nl u\nr_{k,p}=\left(\lVert u\rVert_{p}^{p}+\lVert u'\rVert_p^p+\ldots+\lVert u^{(k)}\rVert_p^p\right)^{1/p}
$$
and denotes the resulting Sobolev space by $W^{k,p}$.
As a first example of what is to come, Peetre recalls an observation by Douady, which shows that the canonical mapping $W^{k,p}\to L^p$ fails to be injective. 
He then proceeds to show that this space has a number of pathological, yet interesting, properties.
The most astonishing is the isomorphism
$$
W^{k,p}\cong L^{p}\oplus L^{p}\oplus \ldots\oplus L^p\cong L^p,
$$
where the last isomorphism is well known. 
This, in turn, gives by a classical theorem of Day that the dual of $W^{k,p}$ is trivial, see \cite{Day}. 
Some remarks concerning this isomorphism are made in \cite{Wennman}, where we in particular discuss the surjectivity of the differentiation mapping, defined on test functions by
$$
(\delta f)(x)=f'(x),\qquad f\in C^\infty(\R).
$$
Many of the results in this paper have analogues for the space $W^{k,p}$.
Naturally, the techniques in \cite{Peetre} and \cite{Wennman} are highly inspiratory for the present work.

\subsec{Carleman-Sobolev classes} 
In \cite{BehmWennman} the author and Gustav Behm investigate several possible extensions of Peetre's investigations. 
Many references both to explicit results and to techniques from that paper are made here, and it is natural to consider this paper as a continuation of it.

Partly inspired by Sobolev spaces and partly by so-called Carleman classes from the study of quasi-analytic functions, we consider what we term Carleman-Sobolev classes; defined as completions of smooth test classes with respect to the weighted norm
\begin{equation}\label{weight-norm}
\nl u\nr_{\M}=\sup_{n\geq0}\frac{\nl u^{(n)}\nr_{p}}{M_n},\quad u\in C^{\infty}(\R)
\end{equation}
where $\M=\{M_n\}$ is a weight sequence. 
We are mostly concerned with a test class which we denote by $\S$ or $\S_\M$, which is defined to consist of all $C^{\infty}(\R)$-smooth functions $u$ with finite $\M$-norm, which furthermore satisfies the growth restriction
\begin{equation}\label{supderiv}
\limsup_{n\to\infty}\nl u^{(n)}\nr_{\infty}^{q^{n}}\leq 1,
\end{equation}
where $q=1-p$ to simplify notation.
We will consider the space resulting from taking the abstract completion of $\S_\M$ with respect to the quasi norm \eqref{weight-norm}. 
This Sobolev-flavoured space is denoted by $\W$. 
In \cite{BehmWennman} it is proved that provided that the weight sequence does not grow too quickly, expressed as
\begin{equation}\label{prod}
\product{M}:=\prod_{n=0}^{\infty}M_n^{q^{n}}<\infty,
\end{equation}
the space $\W$ can be continuously embedded into $C^{\infty}(\R)$. 
Under certain regularity assumptions on the weight sequence, we also demonstrated that this condition is sharp. 
This is to be interpreted in the sense that such an embedding is impossible when $\product{M}=\infty$. 
In fact, even a continuous embedding into $C(\R)$ is impossible.

The conditions \eqref{supderiv} and \eqref{prod} might seem strange. To get an idea where they come from, we refer to the proof of \cite[Proposition~2.2]{BehmWennman}. We remark that this result is due to Hedenmalm, and the techniques are inspired the work \cite{HedenmalmBorichev}.

We shall also work with a space which we denote by $\Winf$, which is the completion of the set of $C^{\infty}$-smooth functions with finite $\M$-norm. 
That is, we do not have any supremum norm control on the derivatives.

We conjectured in \cite[Conjecture~2]{BehmWennman} that the situation is either very controlled, or there is almost a complete lack of structure. 
That is, we can either embed $\W$ into $C^{\infty}(\R)$, or the situation is as bad as what Peetre encounters. 
Loosely speaking, we expect that when $\product{M}=\infty$ we should have a complete independence of derivatives, in the sense that for any $L^{p}$-functions $g_0, \ldots , g_k$ there exists some $f\in \W$ such that $f^{(n)}=g_n$ for $0\leq n\leq k$. 
How to interpret the derivative will be clear from Section~\ref{indep}. 
For the space $\Winf$ we expect such an independence relation to hold for any sequence $\M$ (see e.g. \mbox{\cite[Conjecture~2]{BehmWennman}}). 

The goal of this paper is to resolve these questions, modulo the unexplored matter of what is the right regularity to require of $\M$. 

\newpage

\subsec{Statement of the results} First we recall a few restrictions that we have to put on the sequence $\M$. 
It seems like our tools for studying Carleman-Sobolev classes breaks down when $\M$ behaves too erratically. 
Indeed, the condition~\eqref{prod} is only a summability condition on the sequence. 
This thus allows some derivatives to grow big while strict control are forced upon others, and it is therefore natural to expect that the analysis should be difficult without further assumptions. 
In the case when $\M$ is large, i.e.~$\product{M}=\infty$, the essence of our results is that {\em the space} $\W$ {\em is really big}. 
To substantiate such a claim we aim at approximating $L^{p}$-functions by a method from \cite{BehmWennman}. 
For this method to work we need to require that the logarithmized sequence is convex, that $\log M_n$ grows faster than the second degree polynomial $P(n)=n^2$ and that the sequence $\{q^{n}\log M_n\}$ has a limit. 
This is what we will define as $p$-regular sequences (see Definition~\ref{regular-M}). 
The reasons for why these are to be regarded as regularity assumptions are discussed in more detail in Section~\ref{regularity}. 
Now we can state the main theorem of this paper.
\begin{mainthm}\label{mainintro}
Suppose that $\product{M}=\infty$ and that $\M$ is $p$-regular in the above sense. Then there exists an isomorphism
$$
\W\cong L^p\oplus W^{p}_\MS
$$
where $\MS$ is the shifted sequence $\MS=\{M_{n+1}\}$.
\end{mainthm}
It is immediate that the sequence $\MS$ is $p$-regular and satisfies $\product{M}=\infty$ whenever $\M$ does.
That $\MS$ inherits the properties of $\M$ allows us to iteratively apply Theorem~\ref{mainintro} to extract as many copies of $L^p$ are we would like, resulting in an isomorphism
$$
\W \cong L^{p}\oplus L^p\oplus \ldots L^{p}\oplus W_{\mathcal{M}_n}^p.
$$

This result is to be put in contrast with what happens when $\product{M}<\infty$. 
This case was handled in \cite{BehmWennman} following an observation by Hedenmalm when considering spaces based on functions $f\in C^{\infty}(\R$). Here we modify the proof slightly so as to be able to handle any interval $I$.
\begin{mainthm}\label{embedding}
Suppose that $\product{M}<\infty$. Then there exists a canonical, continuous embedding
$$
\alpha:\W\hookrightarrow  C^{\infty}(I).
$$
\end{mainthm}

In \cite{BehmWennman} we show that at least some growth restriction on the supremum norm of the derivatives in the spirit of \eqref{supderiv} is necessary. 
Indeed, we show that the completion $\Winf$ does not actually depend that much on the weight sequence. 
As long as $\M$ is a sequence of positive numbers one can always embed $L^p$ into $\Winf$. 

From the investigations made in this paper, it is a short step to sharpen the results from \cite{BehmWennman} to obtain the following, which is our second and last main result.
\begin{mainthm}\label{maincor} For any sequence $\M$ of positive numbers, the space $\Winf$ satisfies 
$$
\Winf\cong L^{p}\oplus \Winfshift,
$$ 
where $\MS$ is the shifted sequence $\MS=\{M_{n+1}\}$. 
\end{mainthm}

\subsec{Structure of the paper} In Section~\ref{prel} we will collect material that is to be regarded as folklore and material that has appeared elsewhere, which will be needed in subsequent sections. 
We begin with a discussion of basic properties of $L^p$-spaces for small exponents and continue with the definitions and basic facts about the Sobolev spaces $W^{k,p}$ and the Carleman-Sobolev classes $\W$. 
This is followed by a discussion regarding which weight sequences $\M$ we wish to consider. 
The main method of approximation used here will be described in Section~\ref{spline}.

In Section~\ref{multiplier} we will introduce new tools in the study of Carleman-Sobolev classes that enables us to prove Theorem~\ref{mainintro}. 
The main new ingredient is Lemma~\ref{mainlemma} in Section~\ref{contraction}.

Section~\ref{retraction} is devoted to constructing a retraction $L^{p}\hookrightarrow \W$ of the canonical mapping $\W\to L^p$. 
That is, a bounded linear injection that embeds $L^{p}$ into the space $\W$ which is right-inverse to $\alpha$. 
The work in Section~\ref{retraction} thus puts the state of affairs for $\W$ in agreement with what is already known for $\Winf$ due to results in \cite{BehmWennman}. 
This enables us to use ideas from \cite{Peetre} and \cite{Wennman} to prove our two main results simultaneously in Section~\ref{indep}.

The paper ends with some concluding remarks in Section~\ref{conclusions}, in which we discuss to what extent the problem of characterizing Carleman-Sobolev classes and Sobolev spaces for summability exponents $p$ with $0<p<1$ is solved.

\subsec{Acknowledgements} The author wishes to express his sincerest gratitude to Gustav Behm for many inspiring discussions, and to H{\aa}kan Hedenmalm for very valuable insights and for suggesting this research topic from the beginning. It was H{\aa}kan who realized that the requirements \eqref{supderiv} and \eqref{prod} could ensure an embedding of $\W$ into $C^{\infty}$, and conjectured that violating \eqref{prod} should in some sense bring about an independence of derivatives which we finally prove here, in the setting of regular weight sequences.  

\section{Preliminaries}\label{prel}

\subsection{The spaces} For $p>0$ we consider the space $L^p$ of Lebesgue measurable functions $u$ such that the expression
$$
\nl u\nr_p=\left(\int_{\R}\lvert u(x)\rvert^{p}\dd x\right)^{1/p}
$$
is finite. For any $p$ in this range the space $L^p$ is a complete metric space. 
For $p\geq1$ the expression $\nl \cdot \nr_p$ is a norm, and $L^p$ is a Banach space. 
For those exponents with $0<p<1$, however, the triangle inequality fails. 
What we have is a quasi-triangle inequality
$$
\nl u+v\nr_p\leq K\left(\nl u\nr_p+\nl v\nr_p\right),\quad u,v\in L^p
$$
for some $K>1$ making $L^p$ into a quasi-Banach space. 
Although the usual triangle inequality fails in this setting, a slightly different version survives. What we have at our disposal is the inequality
\begin{equation}\label{triangle}
\nl u+v\nr_p^p\leq \nl u\nr_p^p+\nl v\nr_p^p,\quad u,v\in L^p.
\end{equation}
This we will find especially useful when we wish to apply it iteratively to bound norms of large sums of functions.
 
Let us fix a number $p$ with $0<p<1$. In order to simplify notation we set $q=1-p$. We will consider Sobolev-type spaces for small exponents.

We start off by recalling some definitions of classical Sobolev spaces. In the usual case when $p\geq1$ one often defines the Sobolev space $W^{k,p}(\Omega)$ on a domain $\Omega\subseteq \R^n$ as the collection of all $u\in L^p(\Omega)$ such that all their distributional derivatives up to order $k$ also belong to $L^p(\Omega)$. 
Another way to go at this is to define a space, denoted $H^{k,p}(\Omega)$, as the abstract completion of smooth functions with respect to the Sobolev norm
$$
\nl u\nr_{k,p}=\left(\lVert u\rVert_p^p+ \lVert u'\rVert_p^p+\ldots +\lVert u^{(k)}\rVert_p^p\right)^{1/p}.
$$
These definitions turn out to be equivalent, which is the content of the famous paper \cite{HeqW} with the pithy title $H=W$. 

The multitude of equivalent definitions that one is allowed to choose from in the classical settings disappears when $0<p<1$. 
First of all, one cannot even expect that all definitions generalize properly, and if they do they need not be equivalent anymore.
In \cite{Peetre} Peetre chooses to define the Sobolev spaces $W^{k,p}$ in the latter way. 
That is, as abstract completions of a set of smooth functions with respect to the above Sobolev norm, which is more accurately described as a quasi-norm for this class of exponents.

We shall work with similarly defined spaces, but we considered slightly different, weighted versions of the quasi-norm taking into account derivatives of all orders. 

\begin{definition}\label{CSC}
Assume that $\M$ is a sequence of numbers greater than or equal to one. We define the Carleman-Sobolev class $\W$ as the abstract completion with respect to the norm \eqref{weight-norm} of the test class
$$
\S=\left\{ u\in C^{\infty}(\R): \norm{u}<\infty\and \limsup_{n\to\infty}\nl u^{(n)}\nr_{\infty}^{q^{n}}\leq 1\right\}.
$$ 
\end{definition}

We shall also need the following spaces. They are not the main focus of this paper, but illustrates well the consequences of the growth restriction \eqref{supderiv}.

\begin{definition}\label{CSC}
We define the Carleman-Sobolev class $\Winf$ as the abstract completion with respect to the norm \eqref{weight-norm} of the test class
$$
\C=\left\{ u\in C^{\infty}(\R): \norm{u}<\infty\right\}.
$$ 
\end{definition}

Both these classes obviously depend on the exponent $p$ and the weight sequence $\M$. 

We admit that the growth restriction \eqref{supderiv} might appear to be taken from thin air. However, if one wants things to be at all different from the case of $W^{k,p}$ for finite $k$, we will see that we need at least some such restriction. 
That is, the space $\Winf$ behaves much like Peetre's $W^{k,p}$ no matter which weight sequence $\M$ one considers. 
That it is so follows partly from \cite[Theorem~3.2]{BehmWennman}, and partly from Theorem~\ref{maincor} of this paper. 
To understand what the motivation is for formulating the growth restriction in exactly this way, we advise the reader to study the simple proof of \cite[Proposition~2.1]{BehmWennman}, originally due to H{\aa}kan Hedenmalm.

\subsection{Required regularity of $\M$}\label{regularity} The result \cite[Theorem~3.1]{BehmWennman} serves to provide a converse to \cite[Theorem~2.1]{BehmWennman}, that is, Theorem~\ref{embedding} here. 
The idea is that there should be a phase transition such that $\W$ is very small when \eqref{prod} holds while being very big when \eqref{prod} fails. 
To prove this we needed to put some regularity assumptions on $\M$.

It seems like the phase transition occurs when $\log \M$ is comparable to the sequence $\{q^{-k}\}$. 
Indeed, taking the logarithm of $\product{M}=\infty$ we find that it is equivalent to the statement that the sequence $q^n\log M_n$ is not in $\ell^1$. 
A simple way to ensure that this is the case, is if $q^n\log M_n$ is bounded from below, i.e.~that $\log M_n$ is eventually bigger than $Cq^{-k}$ for some  constant $C>0$. 
In this case, no further regularity is needed.

When $q^n\log M_n$ fails to be in $\ell^{1}$ for some other reason, we will need three regularity assumptions. 
The most basic is logarithmic convexity of $\M$. Moreover, we need convergence of the sequence $\{q^n\log M_n\}$. 
This is actually a bit stronger than what is necessary. What is needed is that $\limsup q^n\log M_n<\infty$ if $\liminf q^n\log M_n= 0$. 
The second regularity assumption says that $\M$ cannot have polynomial growth. 
This seems natural to expect, since if 
$$
\limsup_{n\to\infty} \frac{\log M_n}{P(n)}<\infty
$$
for some polynomial $P$, then $\product{M}<\infty$ automatically follows, contrary to our assumption. 
Therefore we consider it a kind of regularity to assume that also 
$$
\liminf_{n\to\infty} \frac{\log M_n}{P(n)}=\infty.
$$
Actually, we only need the weaker assumption that the above holds for $P(n)=n^2$. We summarize this in the following definition.
\begin{definition}\label{regular-M}
We say that $\M$ is $p$-regular if either
$$
\liminf_{n\to\infty}q^{n}\log M_n>0
$$
or
\begin{enumerate}[(i)]
\item The sequence $\M$ is logarithmically convex,
\item The limit $\lim q^{n}\log M_n$ exists,
\item It holds that \label{polycond}
$$
\liminf_{n\to\infty}\frac{\log M_n}{n^2}=\infty.
$$
\end{enumerate}
\end{definition}

We will need the following result, connected to the third requirement for $p$-regularity.

\begin{proposition}\label{polygrowth} Let $P$ denote a polynomial, and assume that the condition~(\ref{polycond}) above holds. 
Then we have that
$$
\frac{e^{P(n)}}{M_n^{p}}\to 0,\quad n\to\infty.
$$
\end{proposition}
\begin{proof}
Let $\epsilon>0$. 
Next choose a polynomial $Q(n)=An^2$ such that $e^{P(n)}/e^{Q(n)}<\epsilon$ for all $n\geq0$.
By condition (\ref{polycond}) we can find a $n_0>0$ such that whenever $n\geq n_0$ it holds that $\log M_n/n^2 \geq p^{-1}A$, and thus
$\log M_n\geq p^{-1}Q(n)$. It follows that
$$
\frac{e^{P(n)}}{M_n^p}\leq\frac{e^{P(n)}}{e^{p^{-1}Q(n)p}}=\frac{e^{P(n)}}{e^{Q(n)}}<\epsilon,\quad n\geq n_0.
$$
Since $\epsilon$ was arbitrary, the result follows. 
\end{proof}

For future reference we note that $p$-regularity is stable with respect to shifts of the sequence $\M$, since they mainly concern asymptotic behavior. 
If $\M$ satisfies $\product{M}=\infty$ and is $p$-regular, the same will thus hold for $\MS$, and by induction for $\M_i=\{M_{i+n}\}$ for any $i\geq0$.

It is also interesting to note that if $\M$ is $p$-regular and $\product{M}=\infty$, it will also be $r$-regular for $0<r<p$ and satisfy a similar condition as \eqref{prod} where $q=1-p$ is replaced by $s=1-r$. 

\subsection{Canonical mappings} Let $f$ be an element of $\W$. Then $f$ can be represented as a Cauchy sequence $\{f_j\}$ of functions $f_j\in \S_\M$ in the norm $\norm{\cdot}$. That it is Cauchy in this norm implies in particular that $f_j$ is Cauchy in $L^p$, so $f_j$ converges to some function. This function will be unique and we denote it by $\alpha f$. This defines a mapping 
$$
\alpha: \W\to L^p
$$
which is clearly linear and continuous. 

Moreover, we shall need a kind of differentiation mapping, which we will denote by $\delta$. If $\MS$ denotes the shifted weight sequence $\MS=\{M_{n+1}\}$ it is clear that
$$
\sup_{n\geq0}\frac{\nl f_j^{(n)}-f_k^{(n)}\nr_p}{M_{n+1}} \leq \norm{f_j-f_k}\to 0.
$$
It follows that if $W^{p}_{\MS}$ denotes the completion of the test class $\S=\S_\MS$ corresponding to $\MS$ we can define a mapping by
\begin{align*}
\delta:& \W\rightarrow W^{p}_{\MS}\\
	 & \{f_j\} \mapsto \{f_j'\}.
\end{align*}

On the road towards a proof of our main results, we will need to construct a retraction $\beta$ of the canonical mapping $\alpha$. That is, a canonical and continuous injection
$$
\beta:L^p\hookrightarrow \W
$$
satisfying $\alpha\circ\beta=\id$. Moreover we will need it to satisfy $\delta\circ\beta=0$. 

All these mappings appeared in the setting of $W^{k,p}$ in \cite{Peetre}, and the way in which they interplay will turn out to be exactly the same for $W^{k,p}$, $\Winf$ and $\W$ when $\M$ is regular and satisfies $\mu_\M=\infty$.

\subsection{Spline functions as convolutions}\label{spline} 
As mentioned above, we will need to construct a retraction, or right inverse, to the mapping $\alpha$. 
To do this, we shall for each function $f\in L^p$ find a Cauchy sequence of functions in $\S$ with respect to the $\M$-norm such that $f_j\to f$ in $L^p$. 
Moreover we shall need the derivatives $f_j'$ to vanish in $W^p_\MS$.
The elements of $\S$ do not, however, lend themselves to simple computation.
For that reason we shall restrict our attention to a smaller class of functions, about which we have more knowledge.

We will be concerned with a class of functions given by infinite convolution of multiples of characteristic functions on intervals.
These were investigated in previous work, for more details we refer to \cite[Section~3.1]{BehmWennman}.

We consider a decreasing sequence $a=\{a_n\}$ of positive numbers and let
$$
u_{a,n}(x)=H_{a_0}*H_{a_1}*\ldots *H_{a_n}(x), \quad x\in\R
$$
where 
$$
H_{a_n}(x)=\frac{1}{a_n}\chi_{[0,a_n]}(x),\quad x\in\R.
$$
The function $u_{a,n}$ will turn out to be $C^{n}_{0,+}$, that is, $u_n\in C^{n-1}_0$ and the $n$:th derivative, taken in the distributional sense, will be piecewise continuous and compactly supported.
We then set 
\begin{equation}\label{inf-conv}
u_a(x)=\lim_{n\to\infty} u_n(x) = \lim_{n\to\infty}H_{a_0}*H_{a_1}*\ldots *H_{a_n}(x),\quad x\in\R 
\end{equation} 
and obtain a function $u_a\in C^{\infty}_0(\R)$. 
This limiting process will work fine if $a\in\ell^1$, see \cite[Theorem~1.3.5]{Horm1}.

The functions $u_a$ are much easier to handle than the general element of $\S$. 
However, we can do even better. 
In the following section we will recall a result in \cite{BehmWennman} that actually allows us to do many of the computations on the top derivatives $u_{a,n}^{(n)}$ instead of on the derivative $u_a^{(n)}$.
The graph of $u_{a,n}^{(n)}$ has the shape of a finite number of rectangles, and if $\{a_n\}$ satisfies a further condition these rectangles will be disjoint and of equal hight, coinciding with the supremum of $u_a^{(n)}$.
Using this we can estimate $L^p$-norms of $u^{(n)}_a$ effectively in terms of the numbers $\{a_n\}$.
The condition which guarantees this latter fact is the following:
\begin{equation}\label{separation}
\sum_{j>k}a_j\leq (1-c)a_k.
\end{equation}

We want to work with a subset $\mathcal{I}_\M\subset\S_\M$ of function of this type.
To be sure that $u_a$ lies in the class $\S_\M$ we will need to have control on the $\M$-norm and to require that \eqref{limsup} holds. 
In terms of the sequence $a$ the latter is to say that
\begin{equation}\label{limsup}
\limsup_{n\to\infty} \left(\frac{1}{a_0a_1\cdots a_n}\right)^{q^n}\leq 1.
\end{equation}
What we will verify in practice is actually the following which is equivalent to \eqref{limsup}, as is seen by taking logarithms.
\begin{equation}\label{loglimsup}
\limsup_{n\to\infty}q^n\sum_{j=0}^{n}-\log a_j \leq 0.
\end{equation}
We finally define $\mathcal{I}_\M$ as the set
$$
\mathcal{I}_\M=\left\{u_a: \nl u_a\nr_\M<\infty, a\in\ell^1 \text{ is decreasing and satisfies } \eqref{separation} \text{ and } \eqref{limsup}\right\}.
$$
It is an immediate consequence of the definition that $\mathcal{I}_\M\subset \S_\M$.

As already mentioned, these functions were studied in \cite{BehmWennman}. We recall the following proposition from there which will be of utmost importance later on.

\begin{proposition}\label{convol-prop} Let $a=\{a_j\}$ denote a sequence satisfying \eqref{separation}. Then
$$
\nl u_a^{(n)}\nr_{\infty} = \nl u_{a,n}^{(n)}\nr_{\infty}=\frac{1}{a_0a_1\cdots a_n}, \quad n\geq0,
$$
and we have the estimate
\begin{equation}\label{convol-ineq}
\nl u_a^{(n)}(x)\nr_{p}\leq (2-c)^{1/p}\nl u_{a,n}^{(n)}\nr_{p}=(2-c)^{1/p}\frac{2^{n/p}a_n^{1/p}}{a_0\cdot \ldots \cdot a_n},\quad n\geq0,
\end{equation}
where $c$ is the constant $0<c<1$ from \eqref{separation}.
\end{proposition}

We actually have a reversed inequality as well. This will not be used, but it suggests that we do not lose so much when passing to $u_{a,n}$ instead of $u_a$. 
If $c$ again denotes the constant in \eqref{separation} it holds that
$$
\nl u_a^{(n)}\nr_{p} \geq (1-c)^{1/p}\nl u_{a,n}^{(n)}\nr_p.
$$

\section{Multipliers on Infinite convolutions}\label{multiplier}
Our goal is, since time immemorial, to find estimates. 
More specifically estimates for $p$-norms of functions $u_a\in \mathcal{I}_\M$.
Our starting point is the following result from \cite{BehmWennman}.

\begin{proposition}\label{propbw}
If $\M$ is $p$-regular, and $\product{M}=\infty$ then there exists sequences $a_j=\{a_{j,k}\}_{k\geq0}$ such that
$\supp u_{a_j}\to\{0\}$ as $j\to\infty$, $u_{a_j}$ satisfies \eqref{limsup} and
\begin{equation}\label{propbwest}
\nl u_{a_j}^{(n)}\nr_p \leq e^{P(n)}M_n^{q},\quad n\geq0,
\end{equation}
where $P$ is a second degree polynomial with positive leading coefficient.
\end{proposition}
\begin{proof}
It is from this statement that \cite[Theorem~3.3]{BehmWennman} follows.
\end{proof}

\subsection{Existence of mollifiers} The following theorem, guaranteeing the existence of mollifiers with very good convergence properties, will be one of our main tools.
\begin{theorem}\label{mollifier}
Assume that $\M$ is $p$-regular and satisfies $\product{M}=\infty$. For any $\epsilon>0$ there exists a mollifier $v\in\S$ such that
$$
\int_{\R} v\dd x=1, \qquad \supp v\subset [0,\epsilon], \qquad \norm{v}<\epsilon.
$$
\end{theorem}
A very similar result for the class $\C$ instead of $\S$ is found in \cite{BehmWennman}, and was by employing that result a retraction $\beta:L^{p}\hookrightarrow \Winf$ of $\alpha$ could be constructed. 
This suggests that we are in a good position to do the same for the space $\W$. This will be done fairly succinctly in Subsection~\ref{retraction}.
We postpone the proof until the end of the next subsection.

\subsection{Look at another sequence}

We prove Theorem~\ref{mollifier} mainly by using Proposition~\ref{propbw}, as promised. However, applying it directly to the sequence $\M$ will not give us all we need. The following result, however, shows that this fact poses no problem.

\begin{lemma}\label{change}
Let $P$ be any second degree polynomial with positive leading coefficient. If $\M$ is $p$-regular, then there exists a $p$-regular minorant $\mathcal{N}$ to $\M$ such that
\begin{equation}\label{changeest}
\frac{e^{P(n)}N_n^q}{M_n} < \epsilon,\quad n\geq0.
\end{equation}
\end{lemma}

\begin{proof}
The tail is taken care of by Proposition~\ref{polygrowth}, that is for all $n\geq n_0$ for some $n_0$ we have
$$
\frac{e^{P(n)}}{M_n^p}<\epsilon.
$$ 
We can then find finitely many $N_n$ for $0\leq n<n_0$ such that they are log-convex for these $n$, and satisfies
$$
\frac{e^{P(n)}N_n^q}{M_n}<\epsilon.
$$
Next, due to the log-convexity of $\M$ we have enough regularity to be able to find a straight line connecting the point $(n_0, \log N_{n_0})$ to $(n_1, \log M_{n_1})$ for some $n_1> n_0$, in such a way that the slope of $\log M_n$ is greater than the slope of the line at $n_1$. Let $\log N_n$ be defined for $n_0 \leq n < n_1$ by this line, and then let $N_n=M_n$ for $n\geq n_1$. It is clear that $\mathcal{N}$ is log-convex, it satisfies the required estimate for all $n$. Indeed, small and large values are already taken care of, and for $n_0 \leq n < n_1$ we have
$$
\frac{e^{P(n)}N_n^q}{M_n} \leq \frac{e^{P(n)}M_n^q}{M_n}=\frac{e^{P(n)}}{M_n^p}<\epsilon.
$$
It is clear that all other properties required for $p$-regularity of $\mathcal{N}$ are inherited.
\end{proof}

We recall that our goal is to construct functions with integral equal to one, while having very small support and small $\M$-norm. When we apply Proposition~\ref{propbw} directly to the sequence $\M$, we end up with a family of functions whose norms are uniformly bounded, and where the $L^p$-norms of high derivatives are just as we want them. The supports tend to be fine as well. However, using the following proposition we can get all the way.

Now for the proof of our main result in this section.

\begin{proof}[Proof of Theorem~\ref{mollifier}] 
Let $\mathcal{N}$ be the $p$-regular minorant to $\M$ for which the estimate \eqref{changeest} holds for a specified second degree polynomial $P$ (the interested reader can consult the proof of Lemma ~3.4 in \cite{BehmWennman}), guaranteed to exist by Lemma~\ref{change}. Apply next Proposition~\ref{propbw} to the sequence $\mathcal{N}$. We get a family of functions $u_j$ such that the support of $u_j$ tends to $\{0\}$ and such that
$$
\nl u_j^{(n)}\nr_p\leq e^{P(n)}N_n^q.
$$
Let $j$ be big enough for $\supp u_j\subseteq [0,\epsilon]$ to hold, and set $v=u_j$. Then
$$
\nl v\nr_\M = \sup_{n\geq0} \frac{\nl u_j^{(n)}\nr_p}{M_n}\leq\sup_{n\geq0} \frac{e^{P(n)}N_n^q}{M_n}<\epsilon,
$$
by \eqref{changeest}. The function $v$ has integral one, by construction, so the result follows.
\end{proof}

\section{Independence of derivatives}\label{indep}
In this section we aim for a proof of our main results, Theorems~\ref{mainintro} and \ref{maincor}. Theese proofs are alomost identical, so we will find it enough to prove the harder result - Theorem~\ref{mainintro}. The outline will follow that of \cite[Corollary~4.1]{Peetre}. 
The first part can be found in \cite{Peetre}, but there seems to be something missing there (mainly, surjectivity of $\delta$). For the remainder of the proof we refer to \cite{Wennman}. 
The idea is that $\W$ splits as a direct sum $\im \alpha \oplus \im \delta$ since derivatives will be independent, and these mappings are surjective onto their respective co-domains.

\subsection{Construction of a retraction using mollifiers}\label{retraction}
We will begin by studying the canonical mapping $\beta:L^p\to\Winf$, which will be linear, continuous and injective and furthermore will turn out to satisfy
$$
\alpha\circ\beta=\id,\and \delta\circ\beta=0.
$$
The idea is that for any $L^p$-function $f$ we shall find a Cauchy sequence $\{f_j\}\in\W$ such that $f_j\to f$ in $L^p$, but $\{f_j'\}=0$ in $W^p_\MS$.
The method of approximation we will use is convolution by the mollifiers constructed in Theorem~\ref{mollifier}. 
What follows is basically a restatement of Subsection~3.3 from \cite{BehmWennman}, where a similar result was proven for $\Winf$.

We begin with a lemma that immediately takes us much closer to be able to approximate $L^p$ functions with Cauchy-sequences in $\S$.

\begin{lemma}\label{step}
Assume $\product{M}=\infty$ and that $\M$ is $p$-regular. Let $f$ denote a step function
$$
f(x)=\sum_{k\leq N}c_k\chi_{[a_k,b_k]}(x),\quad x\in\R.
$$ 
Then there exists $u=\{u_j\}$ in $\W$, that is a Cauchy-sequence of functions in $\S$, such that $u_j\to f$ as $j\to \infty$ in $L^p$ such that $\delta u=0$.
\end{lemma}

\begin{proof}
Let $\{\epsilon_j\}$ be a sequence decreasing to zero and construct for each $\epsilon_j$ a function $v_j$ by applying Theorem~\ref{mollifier} to the space $W^p_\MS$. Set
$$
u_j(x)=\sum_{k\leq N} c_k v_j*\chi_{[a_k, b_k]}(x),\quad x\in\R.
$$
It is clear that $u=\{u_j\}$ lies in $\W$ and that it has the desire properties. For details, see the proofs of \cite[Lemma~3.8 and 3.9]{BehmWennman}.
\end{proof}
From this it is a matter if verifying some simple properties to obtain the following.

\begin{proposition}\label{beta}
If $\product{M}=\infty$ and $\M$ is $p$-regular, there exists a canonical, linear, continuous and injective retraction \mbox{$\beta:L^p\to\W$} such that
$$
\alpha\circ \beta=\id \and \delta\circ\beta=0.
$$
\end{proposition}

\begin{proof}
For a step function $f$ we construct $\beta f=\{u_j\}$ with the desired properties by the previous proposition. 

If $f\in L^p$ there will exist a sequence $\{f_j\}$ of step functions that converges to $f$. Now, for each such we let $\{u_{j,k}\}_{k\geq0}$ be a representative of $\beta f_j$. Thus
$$
\nl u_{j,k}-f_j\nr_p\to 0\and \nl u_{j,k}'\nr_{\MS} \to 0,\quad k\to\infty.
$$
For each $j$ we let $k=k(j)$ be big enough for both the above quantities to be dominated by $j^{-1}$ and set $v_j=u_{j,k(j)}$. Then let $\beta f=\{v_j\}$.

It is clear that each $v_j$ lies in $\S$, and we have
$$
\nl v_j-v_k\nr_{\M}\leq K\left(\nl v_j-v_k\nr_{p}+\nl v_j'-v_k'\nr_{\MS}\right).
$$
It is $v_j$ tends to zero in $W_{\MS}^p$, and we have
$$
\nl v_j-v_k\nr_p^p\leq \nl v_j-f_j\nr_p^p+\nl f_j-f_k\nr_p^p+\nl v_k-f_k\nr_p^p
$$
which tends to zero by construction of $v_j$ and by the fact that $\{f_j\}$ is Cauchy in $L^p$.

Thus $\beta f$ is an element of $\W$. Moreover, it is clear that $v_j\to f$ in $L^p$. 
What remains to be verified is that the mapping is well-defined, continuous, linear and that it behaves as it should with respect to $\alpha$ and $\delta$.
That it is continuous and linear is evident. 
That it is well-defined is also clear. Indeed, the only choices involved are the choices of the sequences $f_j$ and $u_{j,k}$. 
Had one taken another sequences $\{\tilde{f}_j\}$ and other approximating sequences $\tilde{u}_j$ one would have gotten another sequence $\tilde{v}_j$ for $\beta f$. However, $\{\tilde{v}_j'\}$ would still be zero as an element of $W_\MS^p$ and we would still have
$$
\nl \tilde{v}_j-f\nr_p^p\leq \nl \tilde{v}_j-\tilde{f}_k\nr_p^p+\nl \tilde{f}_j-f\nr_p^p\to0,\quad j\to\infty.
$$
Thus $\{v_j\}$ and $\{\tilde{v}_j\}$ are equivalent.
Lastly we note that by these calculations we immediately get
$$
\alpha \circ \beta f = \alpha \{v_j\}=f
$$
and
$$
\delta \circ \beta f= \delta \{v_j\}=\{v_j'\}=0.
$$ 
Thus the result is proven.
\end{proof}

\subsection{The differentiation mapping $\delta$}
In this section we shall translate ideas from \cite{Peetre} and \cite{Wennman} to prove our main results. 
We have already constructed the essential retraction $\beta$, and will soon be in a good position to do remaining parts of the proof fairly transparently. 
First we need another proposition concerning properties of the mapping $\delta$. 

For each $i\geq0$ we denote by $\M_i$ the shifted sequence $\M_i=\{M_{i+n}\}_{n=0}^{\infty}$. Usually we will drop the index zero and denote $\M_0$ by $\M$, as before.

For this mapping we have the following proposition, which is perfectly analogous to Lemma~3.3 in \cite{Wennman}. 
Since it is perhaps not obvious form the start that the proof carries through, we choose to write it down here as well.

\begin{proposition}\label{delta}
Assume $\product{M}=\infty$ and that $\M$ is $p$-regular. Then $\delta$, seen as a mapping either $\W\to W^p_\MS$ or $\delta:\Winf \to W^{\infty,p}_{\MS}$ is surjective and satisfies
$$
\ker \delta=\im \beta.
$$
\end{proposition}
\begin{proof}
For any given $g\in W_{\MS}^p$ we are to find an element $f\in\W$ such that $\delta f=g$. 
To do this, one is tempted to consider a representative $\{g_j\}$ of $g$, find a primitive $f_j$ to each $g_j$ and hope that $f=\{f_j\}$ will do the trick. 
However, there are a couple of obstacles. 
The first is that if one would define $f_j(x)$ as an integral of $g_j$ over $(-\infty, x)$, then nothing is to stop $f_j$ from being non-zero and constant outside the support of $g_j$. 
Even if that is fixed, so that $f_j$ will surely belong to $L^p$ for each fixed $j$, there is nothing to force $\{f_j\}$ to be Cauchy in $L^p$. 

To overcome the first obstacle we will replace first $g$ with a function $\tilde{g}$ such that $\tilde{g}$ has mean zero but represents the same element in $W^p_\MS$. 
We then define $u_j$ as the primitive of $\tilde{g}_j$ in the above canonical way.
To remedy the second problem of convergence in $L^p$ we use the mapping $\beta$ to lift each $u_j$ to some function $\tilde{u}_j\in\W$ which is very close to $u_j$ in the $L^p$-sense, while having negligible derivatives. 
By letting $\tilde{f}_j=u_j-\tilde{u}_j$ we obtain an element which is mapped to $g$ under $\delta$.

More precicely, we represent $g\in W_\MS^p$ by a Cauchy sequence $\{g_j\}$ of functions in $\S_\M$. 
We then let 
$$
\tilde{g}_j=g_j-\left(\int_{\R}g_j(t)\dd t\right) \cdot \psi_j,\quad j\geq1
$$
where $\psi_j\in \S_{\MS}$ are such that $\int \psi_j=1$ and $\nl g_j\nr_1\nl \psi_j\nr_{\MS}\to 0$ as $j\to\infty$. 
That such $\psi_j$ can be constructed follows from Theorem~\ref{mollifier}.
It is clear that $\tilde{g}_j\in \S_{\MS}$ and that $\int\tilde{g}_j=0$. To see that $\{g_j\}$ represents the same element as $\{g_j\}$ we just observe that
$$
\nl g_j-\tilde{g}_j\nr_{\MS} = \nl\left( \int_{\R}g_j(t)\dd t \right)\psi_j\nr_{\MS}\leq \nl g_j\nr_1 \nl \psi_j\nr_{\MS}
$$
which tends to zero by construction.

Now define a sequence $\{u_j\}$ by setting for each $j\geq1$
$$
u_j(x) = \int_{-\infty}^{x}\tilde{g}_j(t)\dd t,\quad x\in\R.
$$
Since each $u_j\in L^p$, we can consider $\beta u_j\in \W$, where $\beta$ is the retraction from Proposition~\ref{beta}. If $\beta u_j$ has a representative $\{v_{j,n}\}_{n\geq0}$, the $v_{j,n}$ will be elements of $\S_\M$ such that $v_{j,n}$ tends to $u_j$ in $L^p$ while $v_{j,n}'$ tends to zero in $W^p_\MS$ as $n\to\infty$. For each $j$ we can find a large enough $n=n(j)$ such that
$$
\max\left\{\frac{\nl u_j-v_{j,n(j)}\nr_p}{M_0}, \nl v_{j,n(j)}'\nr_{\MS} \right\} \leq j^{-1}.
$$
Now set $\tilde{u}_j=v_{j,n(j)}$ for such an $n(j)$, and define $f=\{f_j\}$ by setting
$$
f_j=u_j-\tilde{u}_j,\quad j\geq1.
$$ 
Recall that we need to check that $f\in\W$ and that $\delta f=g$. To do this we observe that $f_j\in\S$ since it is a linear combination of elements in $\S$, and that
$$
\norm{f_j-f_k}=\max \left\{\frac{\nl f_j-f_k\nr_p}{M_0}, \nl f_j'-f_k'\nr_{\MS}\right\}.
$$
Both quantities within brackets turn to zero. That it holds for the leftmost expression follows from the fact that $\nl f_j\nr_p=\nl u_j-\tilde{u}_j\nr_p/M_0<j^{-1}$. 
That the second expression tends to zero follows since
$$
\nl f_j'-f_k'\nr_\MS\leq K\left( \nl u_j'-u_k'\nr_\MS +\nl\tilde{u}_j'-\tilde{u}_k'\nr_{\MS}\right)\leq K\left(\nl \tilde{g}_j-\tilde{g}_k\nr_\MS + K\max\{j^{-1},k^{-1}\}\right)
$$
and the rightmost expression clearily tends to zero, since $\{\tilde{g}_j\}$ is Cauchy in $W_\MS^p$.

Lastly, take $\{\tilde{g}_j\}$ as a representative of $g$ and observe that
$$
\nl\tilde{g}_j-f_j'\nr_\MS= \nl\tilde{g}_j -u_j'+\tilde{u}_j\nr_\MS=\nl \tilde{g}_j-\tilde{g}_j + \tilde{u}_j'\nr_\MS\to0, \quad j\to\infty.
$$
Thus $\delta f=g$, and the proof is complete for $\W$. 

That it holds for $\Winf$ is proven in exactly the same way. 
The only modifications needed concern which references to use. Instead of Theorem~\ref{mollifier} one will have to refer to the corresponding result for $\Winf$ in \cite{BehmWennman}, which is Lemma~3.7. The result corresponding to Proposition~\ref{beta} is \cite[Theorem~3.6]{BehmWennman}.
\end{proof}

\subsection{Adaptation of Peetre's proof scheme}
\begin{proof}[Proof of Theorems~\ref{mainintro} and \ref{maincor}]
We will present the proof for the space $\W$, that is, we will prove Theorem~\ref{mainintro}. The proof of Theorem~\ref{maincor} is perfectly analogous.

We will prove the theorem by explicitly giving the isomorphism. It will be given by
$$
\phi u=(\alpha u, \delta u),\quad u\in\W.
$$
It is clear that $\phi$ is linear and continuous. 
Recall that we can calculate the norm of $u\in\W$ as
$$
\nl u\nr_\M=\max\left\{ \nl \alpha u\nr_p, \nl \delta u\nr_{\MS}\right\}.
$$
Thus $\phi$ is immediately seen to be injective; if $\varphi u=(0,0)$ then both the quantities $\alpha u$ and $\delta u$ must be zero, so $\nl u\nr_\M=0$.

What remains is to show that it is also surjective. 
To accomplish this we pick an element \mbox{$(g,h)\in L^p\oplus W_\MS^p$} and set out to find an $f \in\W$ such that $\phi f=(g,h)$. 
By the surjectivity of $\delta$ to find an element $f_0\in \W$ such that $\delta f_0=h$. 
We next perturb $f_0$ by an element in $\im \beta$, so that the image under $\delta$ will not change, by setting $f=f_0-\beta(\alpha f_0- g)$. For this function we immediately get
$$
\alpha f = \alpha f_0- \alpha\circ \beta (\alpha f_0-g)=\alpha f_0-\alpha f_0+ g = g
$$
and since $\delta\circ\beta=0$ we get
$$
\delta f= \delta f_0-\delta\circ\beta (\alpha f_0-g)=\delta f_0=h.
$$
Thus $\phi$ is an isomorphism of topological vector spaces, and the proof is complete.
\end{proof}

\section{Concluding remarks}\label{conclusions}

In \cite{Peetre} Jaak Peetre initiated the studies of Sobolev spaces, defined as abstract completions of sets of smooth functions with respect to the Sobolev norm, for small exponents.
Following an observation by Hedenmalm, the author and Gustav Behm proceeded to investigate similar spaces in \cite{BehmWennman}. 
During that work it became clear that both the distinction between the cases $\mu_{\M}<\infty$ and $\mu_{\M}=\infty$, and the growth restriction \eqref{supderiv} are immensely important for the structure of the spaces $\W$ and $\Winf$.

In that paper we conjectured that Theorem~\ref{mainintro} should be true under some regularity assumptions on $\M$. 
That is, there is kind of a phase transition as $\mu_\M$ shifts from finite to infinite. 
At that point we had proven that this condition differentiates between whether it is possible to embed $\W$ into $C^{\infty}(\R)$ or not, but it was still possible that the space $\W$ was not isomorphic to $L^p\oplus W_{\mathcal{N}}^p$ for some other sequence $\mathcal{N}$ with similar properties as our original $\M$, even though $\product{M}=\infty$. Here we show that at least when $\M$ is regular, this isomorphism holds. 

Regarding the regularity, the techniques here do not seem to allow for any substantial improvements. They are needed when constructing the mollifiers $v$ in Theorem~\ref{mollifier}. The underlying construction is not presented here, but in brief what happens is that one constructs a function $v_a$ by infinite convolutions of functions of the type $a_{j}^{-1}\chi_{[0,a_j]}$ for some sequence $a=\{a_j\}$. The function $v_a$ will have small norm if $a$ decays fast to zero, but that will also make the supremum norm of $v_a^{(n)}$ grow large fast as $n$ increases. Thus in the presence of the condition \eqref{supderiv}, one has to get the decay just right. In \cite{BehmWennman} this is done to some extent using the regularity of $\M$ and the assumption $\product{M}=\infty$ very explicitly. 
 Using the machinery from Section~\ref{multiplier} we manage to bend the sequence $a$ into another sequence $b$ so that it is more optimal in the sense that the norm $\nl v_{b}\nr_\M$ is small, while the supremum norms will not start to grow too quickly so as to violate \eqref{limsup}. It's hence the construction of $v_a$ that would have to be different in order for us to be able to drop some regularity assumptions.

The rest of the proof scheme seems like it would survive. That is, if one manages to construct mollifiers with arbitrarily small norm and support one could apply the machinery from Section~\ref{indep} without change, and the multipliers on convolutions will probably still work fine.

Modulo this uncertainty, the author considers the problem presented to be solved. However, there are some issues that are related, which are not yet understood. Here we start out with a class defined by removing some $C_0^\infty$-functions that are in some sense too big, and then we close everything up in our norm. Another way to go at it would be to start with a smaller class for which we know that these requirements are fulfilled. Then we would certainly get a smaller space, but for some classes perhaps this is not strict? For example we could consider the completion of $\mathcal{I}_\M$, since it is those functions we employ to approximate $L^p$-functins by the mapping $\beta$. However, if we instead consider the unit circle $\mathbb{T}$, a more natural example is perhaps the space of trigonometric polynomials. These will certainly belong to $\S_\M$ for any decent $\M$. We thus ask the following question.

\vspace{8pt}

\noindent {\bf Question.} Denote by $TP$ the space of trigonometric polynomials on $[0,2\pi]$, and let $\M$ denote a weight sequence for which $\mu_\M=\infty$. What is the abstract completion of $TP$ with respect to the norm $\norm{\cdot}$?

\bibliography{references}{}
\bibliographystyle{amsalpha}

\end{document}